\documentclass[12pt]{amsart}
\usepackage{amsfonts}
\usepackage{graphicx}
\usepackage{amsmath}
\usepackage{amssymb}
\usepackage{mathtools}
\usepackage{xcolor}
\usepackage{pdfsync}
\makeatletter
\newcommand{\mcite}[1]{\csname b@#1\endcsname}
\makeatother

\theoremstyle{theorem}

 \newtheorem{theorem}{Theorem}
\newtheorem{lemma}{Lemma}[section]
\newtheorem{proposition}{Proposition}
\newtheorem{corollary}{Corollary}

\theoremstyle{definition}
\newtheorem{definition}[lemma]{Definition}

\theoremstyle{remark}
\newtheorem{remark}[lemma]{Remark}
\numberwithin{equation}{section}

\newcommand{\tg}{\mathop{\mathrm{tg}}}
\newcommand{\re}{\mathop{\mathsf{Re}}}
\newcommand{\im}{\mathop{\mathsf{Im}}}
\newcommand{\Smap}{\ell}
\newcommand{\Ham}{{\mathcal{H}}}
\newcommand{\CP}{{\mathbb{C}}}
\newcommand{\UH}{{\mathbb{H}}}
\newcommand{\UD}{{\mathbb{D}}}
\newcommand{\UC}{{\partial\UD}}

\newcommand{\di}{\mathrm{d}}
\newcommand{\DI}{\,\di}
\newcommand{\dd}[2]{\frac{\di #1}{\di #2}}
\newcommand{\ddt}[1]{\dd{#1}{t}}
\newcommand{\Real}{{\mathbb{R}}}
\newcommand{\Hol}{{\sf Hol}}
\newcommand{\Aut}{{\sf Aut}}
\def\id{{\sf id}}
\newcommand{\anglim}{\angle\lim}

\newenvironment{mylist}{\begin{list}{}%
{\labelwidth=2em\leftmargin=\labelwidth\itemsep=.4ex plus.1ex
minus.1ex\topsep=.7ex plus.3ex
minus.2ex}%
\let\itm=\item\def\item[##1]{\itm[{\rm ##1}]}}{\end{list}}

\title[Value regions of univalent self-maps \ldots]{Value regions of univalent self-maps with two boundary fixed points}
\author[P.~Gumenyuk, D.~Prokhorov]{Pavel Gumenyuk${}^\dag$, Dmitri Prokhorov}

\address{P.~Gumenyuk: Department of Mathematics and Natural Sciences, University of Stavanger, N-4036 Stavanger, Norway}
\email{pavel.gumenyuk@uis.no}

\address{D.~Prokhorov: Department of Mathematics and Mechanics, Saratov State University, Saratov 410012,
Russia} \email{ProkhorovDV@info.sgu.ru}

\thanks{${}^\dag$ Partially supported by {\it Ministerio de Econom\'\i{}a y Competitividad} (Spain) project \hbox{MTM2015-63699-P}.}

\begin{document}
\let\ge=\geqslant
\let\le=\leqslant

\newcommand{\UDc}{\hskip.025em\overline{\hskip-.025em\vphantom b\UD\hskip-.1em}\hskip.1em}
\newcommand{\Uone}{{\mathfrak U\hskip.05em}_1\hskip-.07em}
\newcommand{\Gam}[1]{{\displaystyle\Gamma^{\hbox{\mathversion{bold}$\scriptstyle #1$}\!}}}
\newcommand{\Gamp}{\Gam+}
\newcommand{\Gamm}{\Gam{\hskip.07em\raise.2ex\hbox{\mathversion{bold}$\scriptstyle-$}\hskip-.03em}}

\begin{abstract}
In this paper we find the exact value region $\mathcal V(z_0,T)$ of the point evaluation functional $f\mapsto f(z_0)$ over the class of all holomorphic injective self-maps $f:\UD\to\UD$ of the unit disk~$\UD$ having a boundary regular fixed point at~$\sigma=-1$ with $f'(-1)=e^{T}$ and the Denjoy\,--\,Wolff point at $\tau=1$.
\end{abstract}

\maketitle

\section{Introduction}

Since the seminal paper~\cite{CowenPommerenke} by Cowen and Pommerenke, the study of holomorphic functions with finite angular derivative at prescribed boundary points has been an active field of research in Complex Analysis, see, e.g., \cite{MilnVas,BESh,Diagons, Frolova,GoryainovBRFP,PommVas,Vas02}, just to mention some works in the topic.

Given a holomorphic function $f$ in the unit disk $\UD:=\{z:|z|<1\}$ and a point~$\sigma\in\UC$ such that there exists finite angular limit $f(\sigma):=\anglim_{z\to\sigma}f(z)$, the \textsl{angular derivative} at~$\sigma$ is $f'(\sigma):=\anglim_{z\to\sigma}\big(f(z)-f(\sigma)\big)/(z-\sigma)$.

On the one hand, for univalent (i.e., holomorphic and injective) functions~$f$, existence of the angular derivative $f'(\sigma)$ different from~$0$ and~$\infty$ is closely related to the geometry of~$f(\UD)$ near~$f(\sigma)$; moreover, if there exists $f'(\sigma)\neq0,\infty$, then the behaviour of~$f$ at the boundary point~$\sigma$ resembles conformality, see, e.g., \cite[\S\S4.3,\,11.4]{Pommerenke2}.

On the other hand, for the dynamics of a holomorphic (but not necessarily univalent) self-map $f:\UD\to\UD$, a crucial role is played by the points~$\sigma\in\UC$ for which $f(\sigma)=\sigma$ (or,~more generally, $f(\sigma)\in\partial\UD$) and the angular derivative~$f'(\sigma)$ is finite, see, e.g., [\mcite{BRFPLoewTheory}\,--\mcite{ContrMadrigal:AnaFlows}, \mcite{CDP1}, \mcite{CDP2}, \mcite{EGRSh}, \mcite{GoryainovFr}, \mcite{Poggi}]. Such points~$\sigma$ are called \textsl{boundary regular fixed points}, see Section~\ref{S_selfmaps} for precise definitions and some basic theory. In particular, a classical result due to Wolff and Denjoy asserts that if $f\in\Hol(\UD,\UD)$ has no fixed points in~$\UD$, then it possesses the so-called \textsl{(boundary) Denjoy\,--\,Wolff point}, i.e.,  a unique boundary regular fixed point~$\tau$ such that~$f'(\tau)\le1$.

In this paper we study \textit{univalent}  self-maps $f:\UD\to\UD$ with a given boundary regular fixed point~$\sigma\in\UC$ and the Denjoy\,--\,Wolff point~$\tau\in\UC\setminus\{\sigma\}$. Using automorphisms of~$\UD$, we may suppose that $\tau=1$ and~$\sigma=-1$. Our main result is the \textit{sharp value region} of $f\mapsto f(z_0)$ for all such self-maps of~$\UD$ with $f'(-1)$ fixed. To give a detailed statement, fix $z_0\in\UD$, $T>0$ and let $\zeta_0=x_1^0+ix_2^0:=\Smap(z_0)$, where
$$\Smap\colon\UD\to\mathbb S;\quad z\mapsto\log\big((1+z)/(1-z)\big)$$ is a conformal map of~$\UD$ onto the  strip~$\mathbb S:=\big\{\zeta\colon-\pi/2<\im\zeta<\pi/2\big\}$. Define:
$$
a_\pm(T):=\,e^{-T/2}\sin x_2^0\,\pm\,(1-e^{-T/2}),\quad R(a,T):=\log\frac{1-a}{1-a_+(T)}\,\,\log\frac{1+a}{1+a_-(T)},
$$
$$
V(\zeta_0,T):=\Big\{x_1+ix_2\in\mathbb S\colon a_-(T)\le \sin x_2\le a_+(T),~\big|x_1-x_1^0-\tfrac T2\big|\le \sqrt{R(\sin x_2,T)}\,\Big\}.
$$

\begin{theorem}\label{TH_main}
Let $f\in\Hol(\UD,\UD)\setminus\{\id_\UD\}$ and $T>0$. Suppose that:
\begin{mylist}
\item[(i)] $f$ is univalent in~$\UD$;
\item[(ii)] the Denjoy\,--\,Wolff point of~$f$ is $\tau=1$;
\item[(iii)] $\sigma=-1$ is a boundary regular fixed point of~$f$ and $f'(-1)=e^T$.
\end{mylist}
Then
\begin{equation}\label{EQ_belong}
f(z_0)\in\mathcal V(z_0,T):=\Smap^{-1}\big(V(\Smap(z_0),T)\big)\setminus\{z_0\}\quad\text{for any~$z_0\in\UD$}.
\end{equation}

This result is sharp, i.e., for any $w_0\in\mathcal V(z_0,T)$ there exists $f\in\Hol(\UD,\UD)\setminus\{\id_\UD\}$ satisfying \hbox{(i)\,--\,(iii)} and such that ${f(z_0)=w_0}$.
\end{theorem}
We can also characterize functions $f$ delivering boundary points of~$\mathcal V(z_0,T)$. In many extremal problems for univalent functions $f:\UD\to\CP$ normalized by $f(0)=f'(1)-1=0$, the Koebe function $f_0(z):=z/(1-z)^2$ mapping $\UD$ onto $\CP\setminus(-\infty,\tfrac14]$, and its rotations $f_{\theta}(z)=e^{i\theta}f_0(e^{-i\theta}z)$, $\theta\in\Real$, are known to be extremal. For bounded univalent functions $f:\UD\to\UD$ normalized by $f(0)=0$, $f'(0)>0$, the role of the Koebe function is played by the Pick functions~$p_\alpha(z):=f_0^{-1}(\alpha f_0(z))$, $\alpha\in(0,1)$, mapping $\UD$ onto~$\UD\setminus[-1,-r]$, $r=r(\alpha)\in(0,1)$. In our case, it would be natural to expect that some functions of the form $f=h_1\circ p_\alpha\circ h_2$, where $h_1,h_2\in\Aut(\UD)$, are extremal.
\begin{theorem}\label{TH_extremal}
For any $w_0\in\partial\mathcal V(z_0,T)\setminus\{z_0\}$, there exists a unique $f=f_{w_0}$ satisfying  conditions \hbox{(i)\,--\,(iii)} in Theorem~\ref{TH_main} and such that ${f_{w_0}(z_0)=w_0}$. If $w_0=\Smap^{-1}(\zeta_0+T)$, then $f_{w_0}$ is a hyperbolic automorphism of~$\UD$, namely $f_{w_0}(z)=\Smap^{-1}(\Smap(z)+T)$. Otherwise, $f_{w_0}$ is a conformal mapping of $\UD$ onto $\UD$ minus a slit along an analytic Jordan arc~$\gamma$ orthogonal to~$\partial\UD$, with $f_{w_0}'(1)=1$.
Moreover, $f_{w_0}=h_1\circ p_\alpha\circ h_2$ for some $h_1,h_2\in\Aut(\UD)$ and $\alpha\in(0,1)$ if and only if $w_0=\Smap^{-1}\big(x_1^0+\tfrac T2+i\arcsin a_\pm(T)\big).$
\end{theorem}
\begin{remark}
Note that $z_0$ is a boundary point of the value region~$\mathcal V(z_0,T)$, but does not belong to~$\mathcal V(z_0,T)$.
The proof of the above theorem, given in Section~\ref{S_proof-main}, shows that $z_0$ would be included, and this would be the only modification of the value region, if we replaced the equality ${f'(-1)=e^T}$ in condition~(iii) of Theorem~\ref{TH_main} by the \textit{inequality} ${f'(-1)\le e^T}$ and removed the requirement $f\neq\id_\UD$ assuming as a convention that $\id_\UD$ satisfies~(ii). Note also that under the conditions of Theorem~\ref{TH_main} modified in this way, $f(z_0)=z_0$ if and only if~$f=\id_\UD$, see Remark~\ref{RM_no-two-DW}.
\end{remark}

If $f\in\Hol(\UD,\UD)$ has boundary regular fixed points at~$\pm1$, then replacing $f$ by  $h\circ f$, where $h$ is a suitable hyperbolic automorphism with the same boundary fixed points, we may suppose that $\tau=1$ is the Denjoy\,--\,Wolff point. In this way, as a corollary of Theorems~\ref{TH_main} and~\ref{TH_extremal} we easily deduce a sharp estimate for $f'(-1)f'(1)$, which was obtained earlier with the help of the extremal length method in~\cite[Section~4]{Frolova}.
\begin{corollary}
Let $z_0\in\UD$ and let $f\in\Hol(\UD,\UD)$ be a univalent function with boundary regular fixed points at~$1$ and~$-1$. Then
\begin{equation}\label{EQ_estimate}
\sqrt{f'(-1)f'(1)}\ge L\,\big(\sin\im \Smap(z_0),\,\sin \im \Smap(f(z_0))\big),\quad  L(a,b):=\max\left\{\tfrac{1+a}{1+b},\tfrac{1-a}{1-b}\right\}.
\end{equation}
Inequality~\eqref{EQ_estimate} is sharp. The equality can occur only for hyperbolic automorphisms and functions of the form $f=h_1\circ p_\alpha\circ h_2$, $h_1,h_2\in\Aut(\UD)$, $\alpha\in(0,1)$.
\end{corollary}

 Recently, the sharp value regions of $f\mapsto f(z_0)$ have been determined for other classes of univalent self-maps~\cite{JuliaSeb, DVK, OliverSeb}. The main instrument is the classical parametric representation of univalent functions, going back to the seminal work by Loewner~\cite{Loewner}. In this paper, we use a new variant of Loewner's parametric method, which is specific for functions satisfying conditions of Theorem~\ref{TH_main}.  This variant of parametric representation was discovered quite recently, see~\cite{Goryainov-Kudryavtseva, Gumenyuk_parametric}. We discuss it in Section~\ref{S_param}.

 It is also worth mentioning that in~\cite{GoryainovBRFP}, using another specific variant of the parametric representation, Goryainov obtained the sharp value region of $f\mapsto f'(0)$ in the class of all univalent $f\in\Hol(\UD,\UD)$, $f(0)=0$, having a boundary regular fixed point at~$\sigma=1$ with a given value of~$f'(1)$.

To complete the Introduction, we recall another related result announced by Goryai\-nov~\cite{Goryainov-talk}. Dropping the univalence requirement, one can study holomorphic self-maps ${f:\UD\to\UD}$ satisfying conditions~(ii) and~(iii) in Theorem~\ref{TH_main} by using relationships between boundary regular fixed points and the Alexandrov\,--\,Clark measures. In particular, according to~\cite{Goryainov-talk}, the value region $\mathcal D(0,T)$ of $f\mapsto f(0)$ over all such self-maps~$f$ is the closed disk whose diameter is the segment $\big[0,\ell^{-1}(T)\big]$, with the boundary point~$z_0=0$ excluded. Analyzing the functions delivering the boundary points of  $\mathcal D(0,T)$, one can conclude that~$\partial \mathcal D(0,T)~\bigcap~\partial\mathcal V(0,T)=\{0,\ell^{-1}(T)\}$.

\section{Holomorphic self-maps of the unit disk}\label{S_selfmaps}
In this section we  cite some basic theory of holomorphic self-maps of~$\UD$. More details can be found, e.g., in the monograph~\cite{Abate}.

Let $f\in\Hol(\UD,\UD)$ and $\sigma\in\UC$. According to the classical Julia\,--\,Wolff\,--\,Carath\'eodory Theorem, see, e.g., \cite[Theorem~1.2.5, Proposition~1.2.6, Theorem~1.2.7]{Abate}, if
\begin{gather}\label{EQ_dilation}
\alpha_f(\sigma):=\liminf_{\UD\ni z\to\sigma}\frac{1-|f(z)|}{1-|z|}<+\infty,\\ \intertext{then}
\label{EQ_reg-contact}
 \exists~\anglim_{z\to\sigma}f(z)=:f(\sigma)\in\UC,
\quad\exists~\anglim_{z\to\sigma}\frac{f(z)-f(\sigma)}{z-\sigma}=:f'(\sigma)=\alpha_f(\sigma)\frac{f(\sigma)}\sigma,~ \lefteqn{~\text{and}}\\[1.5ex]
\label{EQ_Julia-ineq}
 \frac{|f(z)-f(\sigma)|^2}{1-|f(z)|^2}\le |f'(\sigma)|\,\frac{|z-\sigma|^2}{1-|z|^2}\qquad\text{for all}~z\in\UD,
\end{gather}
with the equality sign if and only if $f\in\Aut(\UD)$. Note that in its turn, existence of the limits in \eqref{EQ_reg-contact} satisfying $f(\sigma)\in\UC$ and $f'(\sigma)\neq\infty$ immediately implies~\eqref{EQ_dilation}.

\begin{definition}\label{DF_contact-fixed}
Points $\sigma\in\UC$ satisfying~\eqref{EQ_reg-contact} are referred to as \textsl{regular contact points} of~$f$. If in addition to~\eqref{EQ_reg-contact},  $f(\sigma)=\sigma$, then $\sigma$ is said to be a \textsl{regular fixed point} of~$f$. The number~$f'(\sigma)$ is called the \textsl{angular derivative} of~$f$ at~$\sigma$.
\end{definition}

Among all fixed points (boundary and internal) of a self-map~$f\neq\id_\UD$, there is one point of special importance for dynamics.
 On the one hand, if $f(\tau)=\tau$ for some $\tau\in\UD$, then by the Schwarz Lemma, $\tau$ is the only fixed point of~$f$ in $\UD$. If in addition, $f$ is not an elliptic automorphism, then  $|f'(\tau)|<1$  and hence  the sequence of iterates $(f^{\circ n})$, $f^{\circ 1}:=f$, $f^{\circ(n+1)}:=f\circ f^{\circ n}$, converges (to the constant function equal) to~$\tau$ locally uniformly in~$\UD$. On the other hand, if $f$ has no fixed points in~$\UD$, then by the  Denjoy\,--\,Wolff Theorem, see, e.g. \cite[Theorem~1.2.14, Corollary~1.2.16,
Theorem~1.3.9]{Abate}, $f$~has a unique boundary regular fixed point $\tau\in\UC$ such that ${f'(\tau)\le1}$ and moreover, $f^{\circ n}\to\tau$ locally uniformly in~$\UD$ as $n\to+\infty$.

\begin{definition}\label{DF_DW-point}
The point~$\tau$ above is referred to as the \textsl{Denjoy\,--\,Wolff point} of~$f$.
\end{definition}

\begin{remark}\label{RM_no-two-DW}
Since the strict inequality holds in~\eqref{EQ_Julia-ineq} unless ${f\in\Aut(\UD)}$,  a self-map~$f$ can have a fixed point in~$\UD$ and a boundary regular fixed point~$\sigma$ with $f'(\sigma)\le1$ only if~${f=\id_\UD}$.
\end{remark}

\begin{figure}[t]
\hbox to\textwidth{\includegraphics[width=0.475\textwidth]{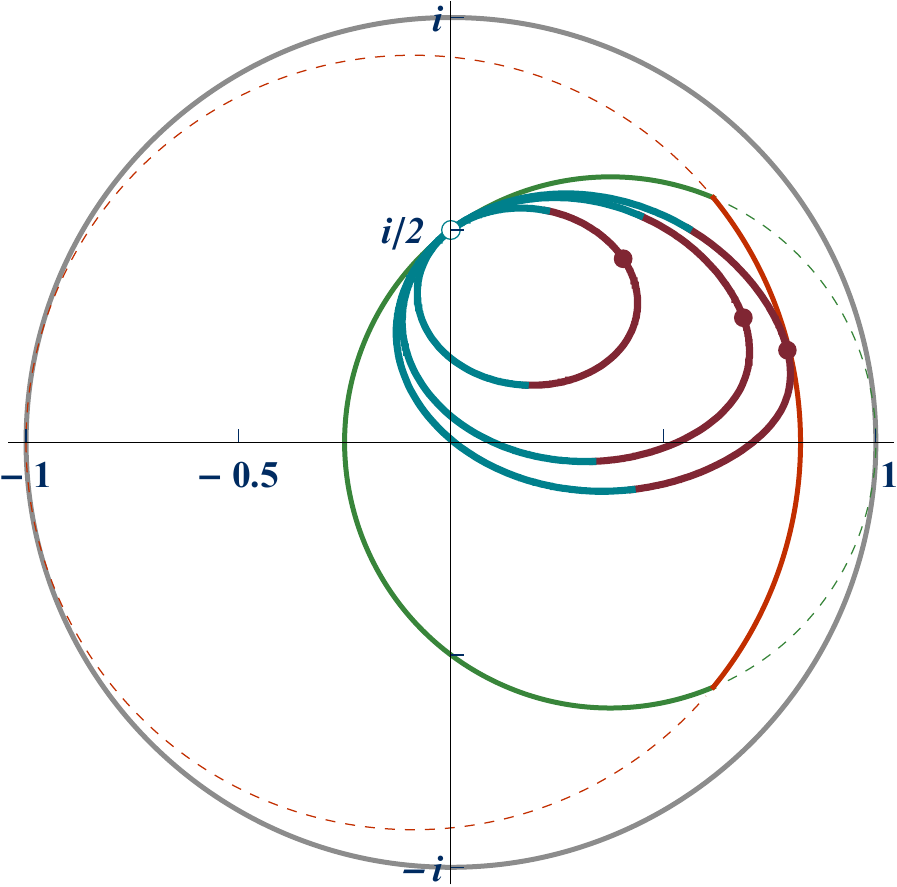}\hss\includegraphics[width=0.475\textwidth]{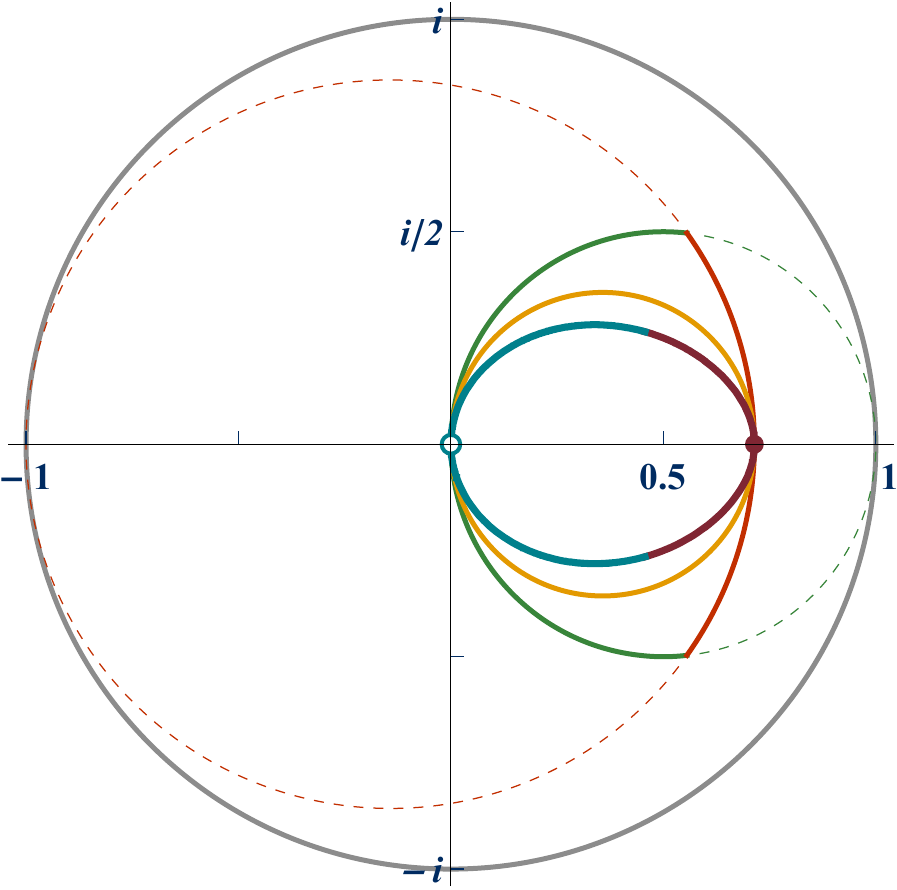}}\vskip.5cm
\definecolor{CGammaM}{rgb}{0, 0.5, 0.55}
\definecolor{CGammaP}{rgb}{0.5, 0.15, 0.2}
\definecolor{Cverde}{rgb}{0.22, 0.52, 0.23}
\definecolor{Croja}{rgb}{0.76, 0.18, 0}
\definecolor{CGoryainov}{rgb}{0.89, 0.6, 0}
\definecolor{CGoryainovM}{rgb}{0.84,0.5,0.01}
\definecolor{UnitDisk}{rgb}{0.55,0.55,0.55}
\begin{center}
\small
\textcolor{UnitDisk}{\raise.5ex\hbox{\rule{.6cm}{1pt}}}~~$\UC$\hskip.8em
 \textcolor{CGammaM}{\hbox to.6cm{\mathversion{bold}\hss$\scriptstyle\circ$\hss}$z_0$}\hskip1.2em
 \textcolor{CGammaM}{\raise.5ex\hbox{\rule{.6cm}{1.25pt}}~~$\Gamm$}\hskip.85em
  \textcolor{CGammaP}{\hbox to.6cm{\mathversion{bold}\hss$\scriptstyle\bullet$\hss}$\ell^{-1}(\ell(z_0)+T)$}\hskip1em  \textcolor{CGammaP}{\raise.5ex\hbox{\rule{.6cm}{1.25pt}}~~$\Gamp$}\hskip1.2em
 \textcolor{Cverde}{\raise.5ex\hbox{\rule{.6cm}{1pt}}~~$\partial D_1$}\hskip1.2em
 \textcolor{Croja}{\raise.5ex\hbox{\rule{.6cm}{1pt}}~~$\partial D_2$}\hskip1.2em
 \textcolor{CGoryainov}{\raise.5ex\hbox{\rule{.6cm}{1pt}}}~~ \textcolor{CGoryainovM}{$\partial \mathcal D(0,T)$}
\end{center}
\caption{{\small The value region~$\mathcal V(z_0,T)$ and the disks~$D_1$, $D_2$ for $z_0:=i/2$, ${T\in\{\log 2,\log 4,\log 6\}}$ and for $z_0:=0$, $T:=\log 6$. The right picture also shows the disk $\mathcal D(0,T)$. Notation $\Gam\pm$ is explained in Section~\ref{S_proof-main}.}}
\label{FG_1}
\end{figure}

\begin{remark}\label{RM_semicont}
Let $f_n(z):=-\sigma H^{-1}\big(\alpha/H(z/\sigma)+\beta/(H(z/\sigma)+n)\big)$, where $n\in\mathbb{N}$,  $\alpha,\beta>0$, and $H(z):=(1+z)/(1-z)$. Note that $f_n(\UD)\subset\UD$ for all~$n\in\mathbb N$ and that $f_n(z)\longrightarrow f(z):=(z+c)/(1+\overline cz)$, $c:=\sigma(1-\alpha)/(1+\alpha)$, locally uniformly in~$\UD$ as $n\to+\infty$. Moreover, $f_n(\sigma)=f(\sigma)=\sigma$ and $f_n'(\sigma)=\alpha+\beta$ for all $n\in\mathbb N$, but $f'(\sigma)=\alpha$. This example shows that the map $f\mapsto f'(\sigma)$ is not continuous.  However, it turns out to be semicontinuous in the following sense. Suppose that $f_n(z)\to f(z)$ as $n\to+\infty$ and that $\sigma\in\UC$ is a boundary regular fixed point of $f_n\in\Hol(\UD,\UD)$ for all~$n\in\mathbb N$ with $\alpha:=\liminf_{n\to+\infty}f_n'(\sigma)<+\infty$. Then passing  in Julia's inequality~\eqref{EQ_Julia-ineq} applied for functions~$f_n$ to the limit, we conclude that $f$ satisfies~\eqref{EQ_Julia-ineq} with $|f'(\sigma)|$ replaced by~$\alpha$. It follows that ${\alpha_f(\sigma)\le \alpha<+\infty}$. Therefore, either $f\equiv \sigma$ or  $f\in\Hol(\UD,\UD)$ and $\sigma$ is a regular boundary fixed point of~$f$ with $f'(\sigma)\le\alpha$. As a consequence, the set of all $f\in\Hol(\UD,\UD)$ sharing two different boundary regular fixed points~$\sigma_1$ and $\sigma_2$ and satisfying $f'(\sigma_j)\le\alpha_j<+\infty$, $j=1,2$, is compact.
\end{remark}

According to inequality~\eqref{EQ_Julia-ineq}, the value region $\mathcal V(z_0,T)$ in Theorem~\ref{TH_main} lies in the intersection of two closed disks $D_1,D_2\subset\UDc$ whose boundaries pass through $z_0$ and $\tau=1$ and through $\ell^{-1}(\ell(z_0)+T)$ and $\sigma=-1$, respectively. Comparison of $\mathcal V(z_0,T)$ with $D_1\cap D_2$ is show in Figure~\ref{FG_1}. On the right picture, for which~$z_0=0$, we also place the value range $\mathcal D(0,T)$ of $f\mapsto f(0)$ over all holomorphic but not necessary injective maps $f:\UD\to\UD$ satisfying conditions~(ii) and~(iii) in Theorem~\ref{TH_main}, see~\cite{Goryainov-talk}.

\section{Parametric representation}\label{S_param}
Denote the class of all $f\in\Hol(\UD,\UD)$ satisfying conditions (i)\,--\,(iii) in Theorem~\ref{TH_main} by~$\mathfrak U(T)$. The following theorem, proved in~\cite{Gumenyuk_parametric}, gives a parametric representation for~$\mathfrak U(T)$ in terms of a Loewner\,--\,Kufarev-type ODE.
\begin{theorem}[\protect{\cite[Corollary~1.2]{Gumenyuk_parametric}}]\label{TH_param}
The class~$\mathfrak U(T)$ coincides with the set of all functions representable in the form $f(z)=w_z(T)$ for all $z\in\UD$, where $w_z(t)$ is the unique solution to the initial value problem
\begin{equation}\label{EQ_LK-specific}
\ddt{w_z}=\tfrac14(1-w_z)^2(1+w_z)q(w_z,t),\quad t\in[0,T],\qquad w_z(0)=z,
\end{equation}
with some function $q:\UD\times[0,T]\to\CP$ satisfying the following conditions:
\begin{itemize}
\item[(i)] for every $z\in\UD$, $q(z,\cdot)$ is measurable on~$[0,T]$;
\item[(ii)] for a.e.~$t\in[0,T]$, $q(\cdot,t)$ has the following integral representation
\begin{equation}\label{EQ_int-repr}
q(z,t)=\int\limits_{\UC\setminus\{1\}}\!\!\frac{1-\kappa}{1+\kappa z}\,\DI\nu_t(\kappa),
\end{equation}
where $\nu_t$ is a probability measure on $\UC\setminus\{1\}$.
\end{itemize}
\end{theorem}
\begin{remark}
A related parametric representation for a class of univalent self-maps of a strip was considered in~\cite{Dubovikov}.
\end{remark}

\begin{remark}\label{RM_param}
In many cases, it is more convenient to deal with the the union $\mathfrak U'(T):=\bigcup_{0\le T'\le T}\mathfrak U(T')$, where we define $\mathfrak U(0):=\{\id_\UD\}$. Indeed, it is evident from the argument of Remark~\ref{RM_semicont} that in contrast to $\mathfrak U(T)$, the class~$\mathfrak U'(T)$ is compact.
Moreover, it is easy to see that Theorem~\ref{TH_param} gives representation of~$\mathfrak U'(T)$ if all probability measures~$\nu_t$ in~\eqref{EQ_int-repr} are replaced with all positive Borel measures~$\nu_t$ satisfying
\begin{equation}\label{EQ_measure}
 \nu_t(\UC\setminus\{1\})\,\in\,[0,1].
\end{equation}
Note that the possibility of $\nu_t=0$ is not excluded.
\end{remark}

\begin{remark}\label{RM_evol-family}
Obviously, the right-hand side of~\eqref{EQ_LK-specific} can be written as $G(w_z,t)$, where
$G(z,t):=\tfrac14(1-z)^2(1+z)q(z,t)$ with $q$ satisfying conditions (i) and~(ii) in Theorem~\ref{TH_param}. By~\cite[Theorem~1]{Goryainov-Kudryavtseva},  $G(\cdot,t)$ is an infinitesimal generator in~$\UD$ for each $t\in[0,T]$.
For simplicity, extend $G$ to all~$t\ge0$ by setting $G(z,t)\equiv0$ for any~$t>T$. Then according to the general theory of Loewner\,--\,Kufarev-type equations, see \cite[Sections\,\,3--5]{BCM1}, for any~$s\ge0$ and any~$z\in\UD$, the initial value problem
$\di w/\di t=G(w,t)$, $t\ge s$, $w(s)=z,$ has a unique solution $w=w_{z,s}(t)$ defined for all $t\ge s$ and the functions $\varphi_{s,t}(z):=w_{z,s}(t)$, $z\in\UD$, $t\ge s\ge0$, form an evolution family, see~\cite[Definition\,\,3.1]{BCM1}.
\end{remark}

\begin{proposition}\label{PR_slit}
Let $\vartheta:[0,T]\to (-\pi,\pi)\setminus\{0\}$, $T>0$, be a $C^1$-smooth function.
 Suppose that in the conditions of Theorem~\ref{TH_param},  $\di\nu_t(e^{i\theta})=\delta\big(\theta-\vartheta(t)\big)\,\di\theta$ for all ${t\in[0,T]}$, where $\delta$ stands for the Dirac delta function. Then $f$ maps~$\UD$ onto $\UD\setminus\gamma$, where $\gamma$ is a slit in~$\UD$, i.e. $\gamma$ is the image of a homeomorphism $\gamma:[0,1]\mapsto\overline\UD$ with $\gamma\big([0,1)\big)\subset\UD$ and $\gamma(1)\in\UC$.

Moreover:
\begin{itemize}
\item[(i)] if $\vartheta$ is a real-analytic function on~$[0,T]$, then $\gamma$ is a real-analytic Jordan arc orthogonal to~$\UD$;
\item[(ii)] $\gamma$ is a circular arc or a straight line segment orthogonal to~$\UC$ if and only if
\begin{equation}\label{EQ_lambdaPick}
\lambda(t):=i\frac{1+e^{i\vartheta(t)}}{1-e^{i\vartheta(t)}}= C_1e^{-t/2}\left(C_2e^{t/2}+\sqrt{C_2^2(e^t-1)+1}\right)^3
\end{equation}
for all~~$t\in[0,T]$ and some constants $C_1,C_2\in\Real$, $C_1\neq0$.
\end{itemize}
\end{proposition}
\begin{proof}
In the conditions of the proposition, \eqref{EQ_LK-specific} takes the following form:
\begin{equation}\label{EQ_LK-specific-slit}
\ddt{w_z}=\tfrac14(1-w_z)^2(1+w_z)\frac{1-e^{i\vartheta(t)}}{1+e^{i\vartheta(t)} w_z},\quad t\in[0,T],\qquad w_z(0)=z.
\end{equation}
The change of variables $\omega_z:=H(w_z)$, where $H(w):=i(1+w)/(1-w)$ maps~$\UD$ conformally onto ${\UH:=\{\omega\colon\im \omega>0\}}$, transforms the above problem to
\begin{equation}
\ddt{\omega_z}=\frac{\omega_z}{1-\lambda(t)\omega_z},\quad t\in[0,T],\qquad \omega_z(0)=H(z),
\end{equation}
where $\lambda(t):=H(e^{i\vartheta(t)})$ for all $t\in[0,T]$. Making further change of variables
$$
\hat\omega_z(t):=\omega_z(t)+\int_0^t\!\frac{\di s}{\lambda(s)},\qquad \xi(t):=\frac{1}{\lambda(t)}+\int_0^t\!\frac{\di s}{\lambda(s)},\qquad \tau=v(t):=\frac{1}{2}\int_0^t\!\frac{\di s}{\lambda(s)^2},
$$
we obtain the chordal Loewner equation
\begin{equation}\label{EQ_chordal}
\frac{\di\hat\omega_z}{\di \tau}=\frac{2}{\xi-\hat\omega_z},\quad \tau\in[0,v(T)],\qquad \hat\omega_z(0)=H(z).
\end{equation}
The geometry of solutions to~\eqref{EQ_chordal} is well-studied, see, e.g.,~\cite{Lind, MarshallRohde, Prokhorov-Vasiliev, IPV,Wong}; see also~\cite{Kufarev46}. In particular, since the function $s\mapsto \xi(v^{-1}(s))$ is $C^1$-smooth, it follows that $z\mapsto \hat\omega_z(T)$ maps~$\UD$ onto $\UH$ minus a slit along some Jordan arc~$\gamma_0$. Taking into account that $w_z(T)=H^{-1}\big(\hat\omega_z(T)-C\big)$,  where $C:=\int_0^T\!\lambda(t)^{-1}\,\di t$, this proves the first part of the proposition.

If $\vartheta$ is real-analytic, then $s\mapsto \xi(v^{-1}(s))$ is real-analytic on~$[0,T]$ as well, and hence by~\cite[Theorem~1.4]{LindTran}, ~$\gamma_0$ is a real-analytic Jordan arc. Moreover, the argument of \cite[Section 6.1]{LindTran} shows that in such a case, $\gamma_0$ is orthogonal to~$\Real$. This proves~(i).

It remains to prove~(ii).  Suppose that $\gamma$ is a circular arc or a straight line segment orthogonal to~$\UC$.
Then we can find a linear-fractional transformation $H_*$ of $\UD$ onto~$\UH$ such that $H_*(\gamma)=[0,i]$. Let~$(\varphi_{s,t})$ be the evolution family associated with equation~\eqref{EQ_LK-specific-slit}, see Remark~\ref{RM_evol-family}.
Note that $\varphi_{t,T}(\UD)\supset\varphi_{t,T}(\varphi_{0,t}(\UD))=\varphi_{0,T}(\UD)=f(\UD)$ for any~$t\in[0,T]$.
It follows that the  intersection of a sufficiently small neighbourhood of $H_*^{-1}(\infty)$ with $\partial\varphi_{t,T}(\UD)$ is an open arc of~$\UC$ containing $H_*^{-1}(\infty)$. Therefore, for each $t\in[0,T]$, there exists a unique $h_t\in\Aut(\UD)$ such that $g_t:=H_*\circ\varphi_{t,T}\circ h_t\circ H_*^{-1}\in\Hol(\UH,\UH)$ satisfies the Laurent expansion $g_t(z)=z-c(t)/z+\ldots$ at~$\infty$ with some~$c(t)\in\Real$.

Denote $H_t:=H_*\circ h_t^{-1}$ for all $t\in[0,T]$. By construction,
$\UH\setminus[0,i]=g_0(\UH)\subset g_t(\UH)\subset g_T(\UH)=\UH$
for all~$t\in[0,T]$.
Thanks to continuity of~$\vartheta$, the function $t\mapsto c(t)$ is $C^1$-smooth.
Therefore, according to the classical result~\cite{Kufarev68} by Kufarev et al, see also~\cite{Andrea}, for any $z\in\UD$, $\tilde\omega_z(t):=g_t^{-1}\circ g_0(H_0(z))$, $t\in[0,T]$, is the unique solution to the initial value problem
$\di{\tilde\omega_z}/\di t=-c'(t)/\tilde\omega_z$, $\tilde\omega_z(0)=H_0(z)\in\UH$.

By construction, $\tilde \omega_z(t)=H_t(w_z)$ for all $t\in[0,T]$ and all~$z\in\UD$. Comparing the differential equations for~$\tilde\omega_z$ and $w_z$, one can conclude that for all $t\in[0,T]$,
\begin{equation}\label{EQ_H_t}
H_t(w):=\frac{\lambda(t)H(w)-1}{a(t)\big(\lambda(t)H(w)-1\big)+b(t)}
\end{equation}
 with real coefficients $a(t)$ and $b(t)$ satisfying
\begin{equation}\label{EQ_sys-a-b}
\di a/\di t=a^3/b^2,\quad \di b/\di t =-3a+b+3a^2/b, \qquad t\in[0,T],
\end{equation}
and such that ${\lambda'(t)/\lambda(t)=1-3a(t)/b(t)}$ and ${b(t)\lambda(t)>0}$ for all $t\in[0,T]$.  System~\eqref{EQ_sys-a-b} can be solved by introducing a new unknown function $k(t):=a(t)/b(t)$. In this way, one can easily check that~$\lambda$ must be of the form~\eqref{EQ_lambdaPick}.

Conversely, if $\lambda$ is given by~\eqref{EQ_lambdaPick}, then system~\eqref{EQ_sys-a-b} has a real-valued solution satisfying ${\lambda'(t)/\lambda(t)=1-3a(t)/b(t)}$ and ${b(t)\lambda(t)>0}$ for all $t\in[0,T]$. It follows that for any $z\in\UD$, the function $\tilde\omega_z(t):=H_t(w_z(t))$, where
$H_t$ is given by~\eqref{EQ_H_t},
is a solution to $\di{\tilde\omega_z}/\di t=-1/\big(b(t)^2\,\tilde\omega_z\big)$, $t\in[0,T]$, $\tilde\omega_z(0)=H_0(z)\in\UH$.
Solving the latter initial value problem for~$\tilde \omega_z$, we conclude that the image of the map ${\UD\ni z\mapsto \tilde\omega_z(T)}$ is the domain $\UH\setminus [0,i\sqrt{Q_T}\,]$, $Q_T:=\smash{2\int_0^T\!b(t)^{-2}\,\di t}$. Thus, $\gamma=H_T^{-1}\big([0,i\sqrt{Q_T}\,]\big)$ is a circular arc or a straight line segment orthogonal to~$\UC$. The proof is now complete.
\end{proof}

\section{Proof of the main results}\label{S_proof-main}
In this section we prove Theorems~\ref{TH_main} and~\ref{TH_extremal}. Fix $T>0$. We start by considering the problem to determine the compact value region $\{f(z_0)\colon f\in\mathfrak U'(T)\}$. Thanks to Theorem~\ref{TH_param} and Remark~\ref{RM_param}, it coincides with the reachable set $\{w_{z_0}(T)\}$ of the controllable system~\eqref{EQ_LK-specific} in which the measure-valued control $t\mapsto \nu_t$ satisfies~\eqref{EQ_measure}.
The change of variables $$\zeta=\Smap(w),\quad\lambda=i\frac{1+\kappa}{1-\kappa},$$
reduces our problem to finding the reachable set~$\Omega'_T:=\{\zeta(T)\}$ for the following controllable system
\begin{equation}\label{EQ_in-half-plane}
\frac{\di\zeta}{\di t}=\int_{\mathbb R}\frac{\di\mu_t(\lambda)}{1-i\lambda e^\zeta},\quad t\in[0,T];\qquad\zeta|_{t=0}=\zeta_0:=\Smap(z_0),
\end{equation}
where $\mu_t$'s are positive Borel measures on~$\Real$ with $\mu_t(\Real)\le1$. By using the prime in the notation $\Omega'_T$ we emphasize that this reachable set corresponds to the class~$\mathfrak U'(T)$.

Denote $x_1:=\re \zeta$ and $x_2:=\im\zeta$. Note that $x_2\in(-\frac{\pi}{2},\frac{\pi}{2})$. For any fixed $\zeta={x_1+ix_2\in\mathbb S}$, the range of the right-hand side in~\eqref{EQ_in-half-plane}, regarded as a function of the measure~$\mu_t$, is the disk $$\left\{\omega\in\CP\colon\left|\omega-\frac{e^{-ix_2}}{2\cos x_2}\right|\le\frac{1}{2\cos x_2}\right\}.$$

Therefore, replacing the measure-valued control $t\mapsto \mu_t$ with the complex-valued control
$$u(t):=2e^{ix_2}\cos x_2\int_{\mathbb R}\frac{\di\mu_t(\lambda)}{1-i\lambda e^{x_1+i x_2}},$$
we can rewrite~\eqref{EQ_in-half-plane} in the following form
\begin{align}
\ddt{x_1}&=\re\frac{u(t)e^{-ix_2}}{2\cos x_2}=\frac12\re u(t)+\frac{\tg x_2}{2}\im u(t), && x_1(0)=x_1^0:=\re\zeta_0, \label{EQ_x1}
\\[1ex]
\ddt{x_2}&=\im\frac{u(t)e^{-ix_2}}{2\cos x_2}=\frac12\im u(t)-\frac{\tg x_2}{2}\re u(t), && x_2(0)=x_2^0:=\re\zeta_0, \label{EQ_x2}
\end{align}
where $u:[0,T]\to U:=\{u\colon|u-1|\le1\}$ is an arbitrary measurable function.

Introduce the Hamilton function
\begin{equation*}
\Ham(x_1,x_2,\Psi_1,\Psi_2,u):=\Psi_1\re\frac{ue^{-ix_2}}{2\cos x_2}+\Psi_2\im\frac{ue^{-ix_2}}{2\cos x_2}=\re\frac{ue^{-ix_2}(\Psi_1-i\Psi_2)}{2\cos x_2},
\end{equation*}
where $\Psi_1$, $\Psi_2$ satisfy the adjoint system of ODEs
\begin{equation}\label{EQ_conj}
\frac{\di\Psi_1}{\di t}=-\frac{\partial \Ham}{\partial x_1}=0,\qquad \frac{\di\Psi_2}{\di t}=-\frac{\partial \Ham}{\partial x_2}=-\im\frac{u(t)(\Psi_1-i\Psi_2)}{2\cos^2 x_2}.
\end{equation}

Boundary points of the reachable set $\Omega'_T$, forming a dense subset of $\partial\Omega'_T$, are generated by the driving functions~$u^*$ satisfying the necessary optimal condition in the form of Pontryagin's maximum principle,
\begin{equation}\label{EQ_max}
\max_{u\in U}\Ham\big(x_1(t),x_2(t),\Psi_1(t),\Psi_2(t),u\big)
         =\Ham\big(x_1(t),x_2(t),\Psi_1(t),\Psi_2(t),u^*(t)\big)
\end{equation}
for all $t\in[0,T]$, see, e.g., \cite{MThOptProc}.
Trajectories $(x_1(t),x_2(t))$ in (\ref{EQ_max}) are optimal in the reachable set problem, and $(\Psi_1(t),\Psi_2(t))$ satisfy the adjoint system~\eqref{EQ_conj} with the optimal trajectories.
In particular, $(\Psi_1(t),\Psi_2(t))$ does not vanish, and hence the maximum in~\eqref{EQ_max} is attained at the unique point $u^*=1+e^{i(x_2+\varphi)}$, where $\varphi:=\arg(\Psi_1+i\Psi_2)$. Therefore, from \hbox{\eqref{EQ_x1}\,--\,\eqref{EQ_conj}} for the optimal trajectories we obtain
\begin{align}
  \label{EQ_x1opt}
 \lefteqn{\ddt{x_1}}\hphantom{\ddt{\Psi_2}}&=\frac{\cos\varphi + \cos x_2}{2\cos x_2}, && x_1(0)=x_1^0,\\
   \label{EQ_x2opt}
 \lefteqn{\ddt{x_2}}\hphantom{\ddt{\Psi_2}}&=\frac{\sin\varphi -\sin x_2}{2\cos x_2}, && x_2(0)=x_2^0,\\
 \ddt{\Psi_1}&=0,\\
  \label{EQ_Psi-opt}
 \ddt{\Psi_2}&=\frac{\sin\varphi-\sin x_2}{2\cos^2 x_2}\,\big|\Psi_1-i\Psi_2\big|.
\end{align}

System~\eqref{EQ_x1opt}\,--\,\eqref{EQ_Psi-opt} is invariant w.r.t. multiplication of $(\Psi_1,\Psi_2)$ by a positive constant. Therefore, we may assume that either $\Psi_1\equiv 0$, or $\Psi_1\equiv 1$, or~$\Psi_1\equiv -1$.

If $\Psi_1\equiv0$, then $\varphi=\pm\pi/2$ and we easily get that for all~$t\ge0$,
\begin{equation}\label{EQ_special-solutions}
x_1(t)=x_1(0)+t/2,\quad \sin x_2(t)=a_\pm(t):=e^{-t/2}\sin x_2(0)\pm(1-e^{-t/2}).
\end{equation}

Now let $\Psi_1\equiv 1$. Then $\varphi\in(-\pi/2,\pi/2)$ and equation~\eqref{EQ_Psi-opt} takes the following form
\begin{equation}\label{EQ_varphi}
\ddt\varphi=\frac{\sin\varphi-\sin x_2}{2\cos^2x_2}\cos\varphi=\frac{\cos\mathrlap{\varphi}\hphantom{x_2}}{\cos x_2}\ddt{x_2}.
\end{equation}
System~\eqref{EQ_x1opt},~\eqref{EQ_varphi} admits the first integral
$$
I(x_2,\varphi):=\frac{1-\sin\varphi}{1+\sin\varphi}\,\,
\frac{1+\sin x_2}{1-\sin x_2}\,>\,0,
$$
and as a result it can be integrated in quadratures. Namely, if $C:=I\big(x_2(0),\varphi(0)\big)\neq1$, we obtain the following identities
\begin{align}\label{EQ_Cx2}
&B_1(t)-CB_2(t)=(C-1)t/2,\\[1ex]
\label{EQ_Cx1x2}
&x_1(t)-x_1(0)=\mathrlap{\frac{B_1(t)-\sqrt CB_2(t)}{\sqrt C-1},}\\[1.5ex]
\nonumber &\quad\qquad\qquad\text{where}\quad B_1(t):=\log\frac{1-\sin x_2(t)}{1-\sin x_2(0)},\quad B_2(t):=\log\frac{1+\sin x_2(t)}{1+\sin x_2(0)}.
\end{align}
Excluding~$C$ from \eqref{EQ_Cx2},\,\eqref{EQ_Cx1x2} and setting $t:=T$ gives
\begin{multline}\label{EQ_x1-of-x2}
x_1(T)=x_1(0)+\frac{1}{2}\Big(T+\sqrt{\big(T+2B_1(T)\big)\big(T+2B_2(T)\big)}\Big)\\=x_1(0)+\frac{T}{2}+\sqrt{R\big(\sin x_2(T),T\big)},
\end{multline}
where we took into account that according to~\eqref{EQ_Cx2}, $$\ddt{}\Big(t+2B_1(t)\Big)=\frac{2C}{1+\sin x_2(t)+C\big(1-\sin x_2(t)\big)}>0$$
and therefore, $T+2B_1(T)>0$.

For $C=1$, we have $\varphi(t)=x_2(t)$ and hence $\di\varphi/\di t=\di x_2/\di t=0$, $\di x_1/\di t=1$. Therefore, if $C=1$, then \eqref{EQ_Cx2} and~\eqref{EQ_x1-of-x2} hold as well.

Since $C>0$, from~\eqref{EQ_Cx2} we obtain that $x_2(T)\in J(T):=\big(\arcsin a_{-}(T),\arcsin a_{+}(T)\big)$. On the other hand, for any $x\in J(T)$ there exists
a unique $C=C(x)>0$ that verifies~\eqref{EQ_Cx2} with $T$ and~$x$ substituted for $t$ and~$x_2(t)$, respectively. Solving $I(x_2(0),\varphi(0))=C(x)$ provides us with the initial condition in equation~\eqref{EQ_varphi} for which $x_2(T)=x$.

Investigating  the case $\Psi_1\equiv-1$ in a similar way, we conclude that $\partial\Omega'_T$ is the union of the two Jordan arcs
$$
 \Gam\pm(T):=\Big\{x_1+ix_2\in\mathbb S\colon a_-(T)\le \sin x_2\le a_+(T),~x_1=x_1^0+\tfrac T2\pm\sqrt{R(\sin x_2,T)}\,\Big\},
$$
which do not intersect except for the common end-points $\omega^\pm:=x_1^0+T/2+i\arcsin a_\pm(T)$, delivered by solutions~\eqref{EQ_special-solutions}.
Taking into account that by the very definition, $\mathfrak U'(T')\subset\mathfrak U'(T)$ for any $T'\in[0,T]$, it follows that $\Omega'_T=V(\zeta_0, T)$.

The next step in the proof is to pass from the class~$\mathfrak U'(T)$ to the class~$\mathfrak U(T)$. In the problem of finding the value region of the functional $f\mapsto f(z_0)$, this is equivalent to replacing the range~$U$ of the admissible controls~$u$ in~\eqref{EQ_x1}\,--\,\eqref{EQ_x2} by $U\setminus\{0\}$. Denote by $\Omega_T$ the corresponding reachable set. By re-scaling the time, the problem to find $\Omega_{T'}$, $T'\in(0,T)$, can be restated as the reachable set problem at the same time~$T$ and for the same controllable system, but with the value range of admissible controls restricted to~$\alpha\big(U\setminus\{0\}\big)$, $\alpha:=T'/T$.
Note also that $\Gamp(T)\cup\Gamm(T)\setminus\{\zeta_0\}\subset\Omega_{T}$ for any $T>0$. Since  $\alpha\big(U\setminus\{0\}\big)\subset U\setminus\{0\}$ for any~$\alpha\in(0,1)$, it follows that
$$
\Gamp(T')\cup\Gamm(T')\setminus\{\zeta_0\}\subset\Omega_{T'}\subset\Omega_T\quad\text{for any~$T\in(0,T]$}.
$$
Thus $\Omega_T=V(\zeta_0,T)\setminus\{\zeta_0\}$, which completes the proof of Theorem~1.

To prove Theorem~2, we have to identify the functions delivering the boundary points of $\mathcal V(z_0,T)$. They correspond to the controls $u^*$ satisfying Pontryagin's maximum principle~\eqref{EQ_max}. It is easy to see from the above argument that every point $\omega\in\partial\Omega'_T\setminus\{\zeta_0\}$ corresponds to a unique control, which is $C^1$-smooth and takes values on~$\partial U\setminus\{0\}$.  It follows that the corresponding measures $\mu_t$ in~\eqref{EQ_in-half-plane} and the measures $\nu_t$ in the Loewner-type representation~\eqref{EQ_LK-specific},\,\eqref{EQ_int-repr} are also unique. They are probability measures concentrated at one point that moves smoothly with~$t$. Namely, $\di\mu_t(\lambda)=\delta\big(\lambda-\lambda^*(t)\big)\,\di\lambda$, where
\begin{equation}\label{EQ_lambda-opt}
\lambda^*(t):=\frac{1-2\cos x_2(t)/\big(e^{-ix_2(t)}+e^{i\varphi(t)}\big)}{ie^{x_1(t)+ix_2(t)}}=e^{-x_1(t)}\, \frac{~\raise.5ex\hbox{$\sin\tfrac{\varphi(t)-x_2(t)}2$}~}{\lower.75ex\hbox{$\cos\tfrac{\varphi(t)+x_2(t)}2$}\,}.
\end{equation}

The point $\omega=\omega_0:=\zeta_0+T\in\Gamp$ corresponds to $C=1$, in which case $\varphi(t)= x_2(t)$ for all~$t\in[0,T]$ and hence $\lambda^*(t)\equiv0$. Therefore, from~\eqref{EQ_in-half-plane} we see that the unique $f\in\mathfrak U(T)$ delivering the boundary point~$\Smap^{-1}(\omega_0)$ of $\mathcal V(z_0,T)$ is the hyperbolic automorphism
$$
 f(z)=\frac{z+c(T)}{1+c(T)z},\quad c(T):=\frac{e^T-1}{e^T+1},\quad\text{for all $z\in\UD$}.
$$

For the common end-points~$\omega^\pm$ of $\Gamp$ and $\Gamm$, which correspond to $\varphi=\pm\pi/2$, formula~\eqref{EQ_lambda-opt} simplifies to
$\lambda^*(t)=\pm e^{-x_1(t)}$. In view of~$\eqref{EQ_special-solutions}$, the latter expression coincides with $\lambda(t)$ given by~\eqref{EQ_lambdaPick}
if we set $C_1:=\pm e^{-x_1^0}$ and $C_2:=0$. Taking into account the correspondence between $\mu_t$ and~$\nu_t$ and applying Proposition~\ref{PR_slit}, we conclude that the unique functions $f\in\mathfrak U(T)$ delivering the points $\Smap^{-1}(\omega^\pm)$ map~$\UD$ onto~$\UD$ minus a slit along a circular arc or a segment of a straight line orthogonal to~$\UC$.

 It remains to compare $\lambda^*(t)$ given by \eqref{EQ_lambda-opt} with $\lambda(t)$ given by~\eqref{EQ_lambdaPick} for the case $\omega\in \partial\Omega_T\setminus\{\zeta_0,\omega_0,\omega^+,\omega^-\}$. Suppose $\omega\in\Gamp\setminus\{\omega_0,\omega^+,\omega^-\}$. Using equations \eqref{EQ_x1opt},\,\eqref{EQ_x2opt},\,\eqref{EQ_varphi}
and taking into account the first integral $I(x_2,\varphi)=C$, we find that
$$
 \Big(\,1\,+\,2\,\frac{\di}{\di t}\log\lambda^*(t)\,\Big)^2= \left(\frac{\,\cos\lefteqn{\varphi(t)}\hphantom{x_2(t)}}{\,\cos x_2(t)}\right)^2=\frac{C\big(1-a^2\big)}{\big((1+C)a+(1-C)a^2\big)^2},\quad a:=\sin x_2(t),
$$
while $\big(1+2({\di}/{\di t})\log\lambda(t)\big)^2=9C_2^2\,e^t/\big(1+C_2^2(e^t-1)\big)$. However, according to~\eqref{EQ_Cx2}, $e^t$~cannot be expressed as a rational function of~$\sin x_2(t)$. This shows that $\lambda^*$ is not of the form~\eqref{EQ_lambdaPick} and hence, by Proposition~\ref{PR_slit}, the unique function $f\in\mathfrak U(T)$ that delivers the boundary point $\Smap^{-1}(\omega)$ maps $\UD$ onto $\UD$ minus a slit along a real-analytic arc~$\gamma$ orthogonal to~$\UC$ but different from a circular arc or a segment of a straight line. A similar argument applied to the case $\omega\in\Gamm\setminus\{\zeta_0,\omega^+,\omega^-\}$ completes the proof of Theorem~\ref{TH_extremal}.\qed

\end{document}